\documentclass{elsarticle}

\usepackage[utf8]{inputenc}
\usepackage[]{todonotes}
\usepackage[ruled]{algorithm2e}
\usepackage{endnotes}
\usepackage{pgfplots}
\usepackage{lscape}
\usepackage{setspace}
\usepackage{url}
\usepackage{adjustbox}
\pgfplotsset{compat=newest}

\usepackage{verbatim}

\usepackage{amsmath,amsfonts,amsthm,amssymb, mathrsfs}

\usepackage{footnote}
\makesavenoteenv{tabular}

\graphicspath{{./figures/}}

\usepackage{comment}
\usepackage{graphics}
\usepackage{graphicx}
\usepackage{adjustbox}
\usepackage{tikz}
\usepgfplotslibrary{dateplot}
\usepgfplotslibrary{groupplots}
\usepackage{caption}
\usepackage{subcaption}
\usepackage{longtable}
\usepackage{psfrag}
\usepackage{soul}
\usepackage{mwe}

\usepgfplotslibrary{fillbetween}

\usetikzlibrary{er,positioning, calc, shapes, arrows, chains, fit}

\definecolor{myblue}{RGB}{19,145,215}
\definecolor{mygreen}{RGB}{80,176,50}
\definecolor{myred}{RGB}{222,36,16}

\usepackage{lipsum}
\makeatletter
\def\ps@pprintTitle{%
 \let\@oddhead\@empty
 \let\@evenhead\@empty
 \def\@oddfoot{}%
 \let\@evenfoot\@oddfoot}
\makeatother

\begin{document}

  \begin{frontmatter}

\title{Solving the vehicle routing problem with deep reinforcement learning} 

\author{Simone Fo\`{a}$^{\text{a},1}$\footnote{$^1$ Master of Science in Management engineering. }}
\ead{foa.1803733@studenti.uniroma1.it}

\author[1]{Corrado Coppola}
\ead{corrado.coppola@uniroma1.it}

\author[2]{Giorgio Grani}
\ead{g.grani@uniroma1.it}

\author[1]{Laura Palagi}
\ead{laura.palagi@uniroma1.it}

\address[1]{Sapienza University of Rome, Dep. of Computer Science, Control and Management Engineering, Rome, Italy
}
\address[2]{Sapienza University of Rome, Dep. of Statistical Science, Rome, Italy
}

\begin{abstract}
Recently, the applications of the methodologies of Reinforcement Learning (RL) to NP-Hard Combinatorial optimization problems has become a popular topic. This is essentially  due to the nature of the traditional combinatorial algorithms, often based on a trial-and-error process. RL aims at automating this process. At this regard, this paper focuses on the application of RL for the Vehicle Routing Problem (VRP), a famous combinatorial problem that belongs to the class of NP-Hard problems.

In this work, as first, the problem is modeled as a Markov Decision Process (MDP) and then the PPO method (which belongs to the Actor-Critic class of Reinforcement learning methods) is applied. In a second phase, the neural architecture behind the Actor and Critic has been established, choosing to adopt a neural architecture based on the Convolutional neural networks, both for the Actor and the Critic. This choice resulted effective to address problems of different sizes.

Experiments performed on a wide range of instances show that the algorithm has good generalization capabilities and can reach good solutions in a short time. Comparisons between the algorithm proposed and the state-of-the-art solver $OR-TOOLS$ show that the latter still outperforms the Reinforcement learning algorithm. However, there are future research perspectives, that aim to upgrade the current performance of the algorithm proposed.

\end{abstract}
\end{frontmatter}
\parindent 0cm

\section{ Introduction}
\label{sec:Introduction}
The Vehicle Routing problem (VRP) is among the most studied combinatorial Optimization problems due to its application interest. In its easiest version it consists, given a set of nodes and a depot, in determining the set of routes at the minimum cost, where each route must start and end at the depot and all nodes must be visited. The problem was introduced by Dantzig, George B. and Ramser, John H. (\cite{dantzig1959truck}) in 1959  to address the problem of the optimal routing of a fleet of gasoline delivery trucks between a bulk terminal and many service stations. Through the time many versions of the problems have been introduced (see \cite{braekers2016vehicle}) and nowadays the so called green vehicle routing problem (see \cite{lin2014survey}) is gaining particular importance because of the increased attention to environmental issues.

The VRP can be seen as an extension of the Travelling Salesman Problem (TSP \cite{junger1995traveling}) and for this reason, it is proved to be NP-complete (\cite{papadimitriou1977complexity}). Hence, the numerous attempts to solve the problem exactly through mathematical formulations and using Branch and X methods (branch and bound, branch and cut, branch and price) do not provide effective results when the  dimension of the instances becomes high. For this reason, much of the research has shifted toward the development of heuristics (See paragraph 2).

Recently a new approach for combinatorial problems, based on the application of deep reinforcement learning, is gaining a remarkable attention. This is due to the ability of RL to learn through trials and errors (and the trial-and-error process is indeed a typical part of the standard heuristics). The paper \cite{khalil2017learning} represents a milestone for the application of RL for combinatorial optimization problems. In this paper, the authors apply a constructive heuristic for three different combinatorial problems, showing the effectiveness of the solution proposed across multiple problem configurations.
Regarding the application of reinforcement learning framework for VRP, many works have been proposed. The paper of Nazari \cite{nazari2018reinforcement} represents a landmark: in this work, the authors proposed to tackle the problem using a constructive heuristic based on reinforcement learning and they approximate the policy using a recurrent neural network (RNN) decoder coupled with an attention mechanism. At this regard, the work \cite{kool2018attention} proposed an approach to tackle the problem based on attention layers as well. Lastly, the paper \cite{lu2019learning} is among the few that proposed an improving reinforcement learning heuristic rather than a constructive one.

In this context, in our work we developed an improving heuristic based on deep reinforcement learning, which differs from all others because it tackles the problem in two stages (the first stage is the assignment and the second is the routing, inspired by the classic heuristic approach Cluster first- Route second). The RL algorithm proposed aims at learning just an assignment strategy, while the routing part is performed using a traditional heuristic that guarantees a good estimate of the optimal solution. The RL algorithm proposed belongs to the class of the Actor-Critic methods and the neural architecture used to approximate the policy and the value function is constituted by Convolutional neural networks (CNN,for a detailed description see \cite{lecun2015deep}).

This paper is organized as follows: after a brief introduction regarding the Vehicle routing problem and the basics of Reinforcement learning in section (2), the section (3) is dedicated to the description of how we modeled the VRP as a Markov Decision Process. Then, in section (4), we describe in detail the deep neural networks used. The section (5) is dedicated to the description of the algorithm, showing the pseudo-code. In section (6) we discuss the computational experience and finally in section (7) we report the conclusions.


\section{Preliminaries and notation}\label{sec:preliminaries_and_notation}
In this section, we describe the Vehicle routing problem from a mathematical optimization point of view and the basics of reinforcement learning (RL), focusing on one class of methods (The actor critic methods). Among this class of methods, the Proximal Policy Optimization algorithm (PPO) is discussed.
\subsection{The vehicle routing problem}
The vehicle routing problem is a famous combinatorial optimization problem,that has a tremendous application interest as it is faced on a daily basis by thousands of distributors
worldwide and has significant economic relevance. There are many versions of this problem, that try to adapt the optimization problem to a particular application. However, all the versions are extensions of the classic VRP. Among all the versions, the most researched, according to the literature, is the CVRP (Capacitated vehicle routing problem), which is the main focus of this section. 
\begin{figure}[h]
    \centering
    \includegraphics[scale=0.4]{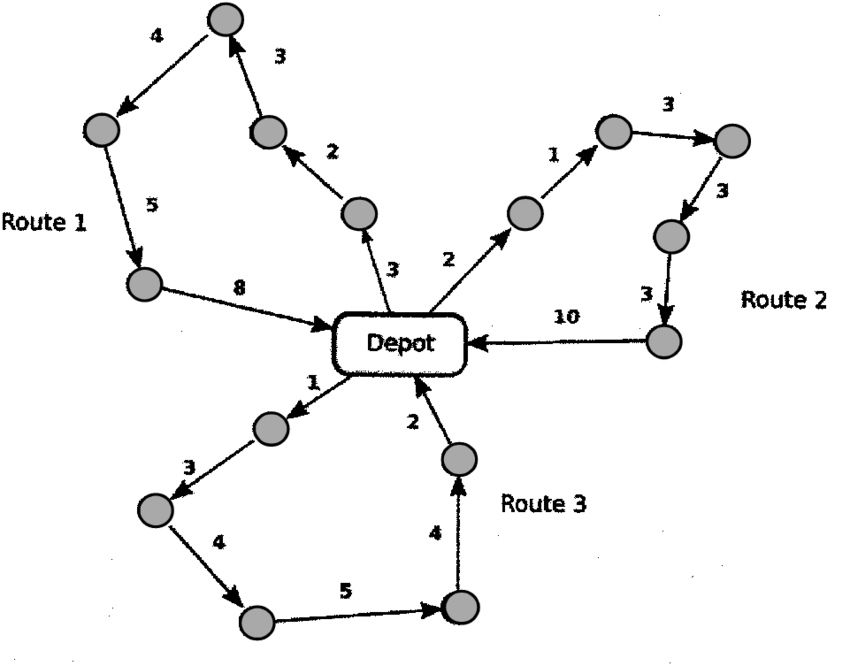}
    \caption{An example of VRP with $n=16,m=3$}
   \label{fig: VRP}
\end{figure}
Let $G=(N,A)$ be a complete graph, where $N=\{0,1,...,n\}$ is the set of nodes, while $A = \{(i,j): i,j \in N,i\neq j\}$ is the set of arcs. The node 0 represents the depot. Let $V=\{1,2,...,m \}$ be the set of vehicles. In addition,let $F=\{1,2,...,l\}$ be the set of features. Each node $i \in N\setminus 0$ is characterized by a non-negative demand of a feature $f\in F$ $d_{if}$. Each vehicle $v \in V$ is characterized by a specific capacity with respect to a particular feature $f\in F$  $q_{vf}$. Moreover, a cost is associated to each arch of the graph $c_{ij}$, $i,j\in N$. In this specific setting, a symmetric cost structure is assumed, so that $c_{ij}=c_{ji} \forall i,j \in N$. The problem consists of determining the sets of $m$ routes with the following properties: (1) each route must start and end at the depot, (2) each node is visited exactly by one vehicle, (3) the sum of the nodes demands for each feature, given a single route, cannot exceed the capacity of the vehicle for that feature, (4) the total routing cost is minimized. The VRP problem is NP-hard because it includes the Traveling salesman problem (TSP) as a special case when the number of vehicles $m=1$ and when $q_{v_f}=\infty,\forall v\in V,f\in F$.

Several families of exact algorithms were introduced over time for the VRP. These are based, first, on the formulation of the problem through a linear programming formulation and then they typically apply either the Branch and X methods (Branch and bound, Branch and price, Branch and cut) or the Dynamic programming. According to in \cite{laporte2007you},the main formulations are the following: (1)Two-Index-Vehicle flow formulation,which has the disadvantage of having a number of constraints exponential with respects to the number of nodes. For this reason, the exponential constraints are generated dynamically during the solution process as they are found to be violated. The second formulation is (2) the set partitioning, which is computationally impractical due to the large number of variables. For this reason, Column generation is a natural methodology applied for this type of formulation. Due to the prohibitive computational effort when the number of nodes becomes large, much of the research regarding VRP focused on the development of heuristics rather than exact methods.

Among heuristics, the  main distinction is between classical heuristics and metaheuristics. The former, at each step, proceed from a solution to a better one in its neighborhood until no further gain is possible. On the other hand,metaheuristics allow the consideration of non-improving and even infeasible intermediate solutions. Regarding the classic heuristics, the most popular is the Clarke and Wright \cite{clarke1964scheduling} (which is an example of constructive heuristics i.e. the algorithm starts form an empty solution and construct it step by step). Its popularity is due to the speed and ease of implementation rather than its accuracy. Among the classic heuristics, an important class of methods is constituted by the improving heuristics, which start from a feasible solution and seek to improve it.
Regarding metaheuristics, they can be classified in (1)local search, (2) population search and (3) learning mechanism. The following work falls into the latter category.

\subsection{The reinforcement learning}
Reinforcement learning (RL) is a paradigm of machine learning, along with supervised learning and unsupervised learning. With regard to the other two paradigms, RL is the most comparable to human learning due to the fact that the new knowledge is acquired  dynamically through  a trial-and-error process. In order to describe effectively the framework of reinforcement learning, the following elements must be described: the agent, the environment, the state, the action, the reward. The agent is the decision maker, who is integrated into an environment. The state represents the encoding of the position in which the agent is located. Given a state, the agent can perform several actions, which cause it to change state. The combination of state and action results in a reward for the agent. In this framework, the goal of the agent is to take suitable actions in order to maximize the cumulative reward. When a model for the environment is present (which means the rules that define the dynamics of states and state-action-reward combinations are known), RL can be formalized as a Markov Decision Process. In this context, the goal of maximizing the cumulative reward can be performed, in theory, using the tools of dynamic programming (\cite{bellman1966dynamic}). However, in real world applications, the full tree of possibilities an agent may encounter cannot be inspected since the size of the tree would be too large to make the computation feasible. Hence, the goal of RL becomes to maximize the expected total reward, so that this can be obtained by inference by sampling a portion of the decision space.

From a formal point of view, we define an episode as an instance to be solved by a RL algorithm. The episode is composed by $T+1$ steps, indexed by $t=0,1,...T$. To show the dependency upon the steps, we use the notation $s_t$ for the states, $a_t$ for the actions, $r_t$ for the reward. In addition, we define $\emph{S}$ the set of all possible states, $\emph{A}$ the set of all possible actions and $\emph{R}:\emph{S}\times \emph{A}\rightarrow \mathbb{R}$ is a function that maps from states and actions to rewards. Furthermore, the policy is identified by the function $\pi(a_t|s_t)$, which represents the probability of taking the action $a_t$ when the state $s_t$ is observed.

Then,the objective function to be maximized is the following:

\begin{equation}
    \label{eq_1}
    J = \mathbb E_\pi \left[ \sum_{t=1}^{T} \gamma^t r_{t} \right]
\end{equation}

$\mathbb E_\pi$ is the expectation of the cumulative discounted reward according to the policy distribution $\pi$ and $\gamma\in (0,1)$ is the discounted factor for the future rewards. The formulation of the objective function \eqref{eq_1} doesn't suggest an intuitive algorithmic approach to solve the problem, therefore many auxiliary functions have been introduced to quantify the consequences of the actions taken at a given step.

\begin{figure}[h]
    \centering
    \includegraphics[scale=0.5]{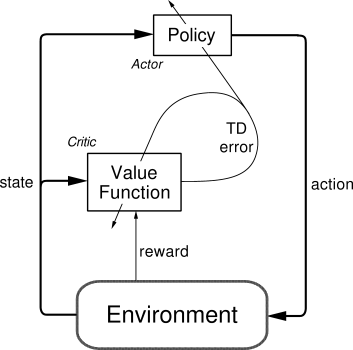}
    \caption{Actor Critic framework}
    \label{fig: The actor critic framework}
\end{figure}

\begin{itemize}
    \item The state-value function. This function represents goodness of a given state $s_t$ (expressed in terms of expected discounted rewards), according to the current policy $\pi$. It can be written as follows.

    \begin{equation}
    \label{eq:Value_Function}
        V_{\pi} (s_t) = \mathbb E_{\pi} \left[\sum_{k=0}^{T-t} \gamma^k r_{t+k}  | S_t = s_t \right] \ \ 
    \end{equation}

\item The action-value function. This function represents the quality of taking a certain action $a_t$ in a given state $s_t$. Even in this case the quality is expressed as the expected discounted reward, according to the current policy $\pi$. It can be written as follows.
\begin{equation}
    \label{eq:Quality func}
        Q_{\pi} (s_t,a_t) = \mathbb E_{\pi} \left[\sum_{k=0}^{T-t} \gamma^k r_{t+k}  | S_t = s_t, A_t = a_t \right]
\end{equation}
\item The advantage function. It is a function that  measures how much is a certain action $a_t$ a good or bad a certain state $s_t$. It returns values greater than $0$ if $Q_{\pi}(s_t,a_t)$ is greater than  $V_{\pi} (s_t)$ It can be written as follows.
\begin{equation}
    \label{eq:AdvF}
        A(s_t,a_t) = Q_{\pi}(s_t,a_t) - V_{\pi} (s_t)
    \end{equation}
\end{itemize}
Over the course of time, several types of RL algorithms
have been introduced and, according to \cite{grondman2012survey},  they can be divided into three
groups: actor-only, critic-only and actor critic methods. The actor-only methods typically parameterize the policy so that classical optimization procedures can be performed. The drawback of this class of methods is the high variance in the estimates of the gradient, leading to slow learning. Among this class, the main two methods are REINFORCE and Policy Gradient Methods. On the other hand, Critic-only methods, such as Q-learning and
SARSA , approximates the state-action value function and no explicit function for the policy approximation is present (See \cite{sutton2018reinforcement}). The actor-critic methods, instead,are characterized by having the agent separated into two decision entities: the
actor and the critic. The critic approximates the state-value function $\hat{V}(s)$, while the actor improves the estimate of the stochastic policy $\hat\pi$ by taking into account the critic estimation. The Actor and Critic do not share weights: therefore, we indicate $\theta$ the weights of the Actor Model and $\omega$ the weights of the Critic model.
In order to optimize  the eq \ref{eq_1}, Actor-Critic methods use the Policy Gradient theorem \cite{sutton2018reinforcement}.The theorem, expressed in the form of log probabilities, states that:
\begin{equation}
    \label{policy_gradient_theorem}
    \begin{aligned}
    \nabla_\theta J_\theta \propto 
    \mathbb{E}_{\pi}\left [ \sum_{a_t} \log \nabla_ \pi(a_t|s_t) A_{\omega} (a_t,s_t)\right]
    \end{aligned}
\end{equation}
This result shows how the gradient depends either by the Actor or the Critic and suggests an algorithmic approach for solving the problem based on gradient ascent. An abstract representation of a generic Actor-Critic algorithm is presented in the figure \ref{fig: The actor critic framework}.

The actor-critic approach used in this work is an adaptation of the Proximal Policy Optimization (PPO) algorithm, more specifically the Adaptive KL Penalty Coefficient version \cite{schulman2017proximal}.
According to the method, the actor updates its parameters through the maximization of the following objective function:

\begin{equation}
    \label{Actor objective function}
    \begin{aligned}
    \max_\theta
    \mathbb{E}_{\pi_\theta}\left[  \frac{\pi_\theta({a_t} | {s_t})}{\pi_{\theta_{old}}(a_t | s_t)} A_{\omega} (a_t , s_t)- \beta KL(\pi(.| s_t),\pi_{old}(.| s_t))\right]
    \end{aligned}
\end{equation}
and the parameter $\beta$ is selected  as follows:
\begin{itemize}
    \item Compute $d=\mathbb{E}[KL[\pi(.|s_t),\pi_{old}(.|s_t]]$
    \item Initialize $d_{targ}$ heuristically
    \item Initialize $\beta$ heuristically
    \item If $d< \frac{d_{targ}}{1.5}, \beta \leftarrow \frac{\beta}{2}$; \text{                         }                  
    If $d> d_{targ}\cdot 1.5, \beta \leftarrow 2\cdot \beta$

\end{itemize}
The parameters of the Critic, instead, are updated minimizing a mean square error loss function, which is the following:
\begin{equation}
    \label{Critic objective function}
    \begin{aligned}
    \min_\omega
    \mathbb{E}_t\left[V_\omega(s_t)-V_{target}\right]^2
    \end{aligned}
\end{equation}
where $V_{target}$ is the sum of the discounted rewards collected.
Hence, a PPO algorithm works as follows: for each iteration, $N$ parallel roll-outs composed by $T$ steps are performed,collecting $N\cdot T$ samples. These are used to update the Actor parameters $\theta$ according to (\ref{Actor objective function}) and to update Critic parameters $\omega$ according to (\ref{Critic objective function}) .

\begin{algorithm}\caption{Proximal Policy Optimization}\label{ref:PPO_Algorithm}
\DontPrintSemicolon
\SetAlgoLined
 Initialize $\theta=\theta_0$, $\omega=\omega_0$, $\beta=\beta_0$;
 {

 \For{iteration =1, \dots, $k$}{
 \For{rollout = 1, \dots, N}{
  Run policy  $\pi_{old}$ in environment for $T$ timesteps \;
  Compute the advantage function estimate $\hat A_1 \dots \hat A_T$\;
 }
 {
 Update the Actor parameter $\theta$  according to (\ref{Actor objective function})  \;
 Update the Critic parameter $\omega$  according to (\ref{Critic objective function})  \;

 }
 
 }

 }
\end{algorithm}

\section{Vehicle routing problem as a Markov decision Process}
\label{sec:Vehicle routing problem as a Markov decision Process}
In this section, we describe the approach used to model the vehicle routing problem as a Markov Decision Process (MDP). VRP and combinatorial problems in general, as already mentioned in the paragraph 2, have been historically addressed using exact algorithms and heuristics.

The approach used in this paper, despite being based on reinforcement learning, is inspired by two principles typical of the traditional heuristics: the improving mechanism and the two stage method. Indeed the VRP can be divided in two stage: the assignment and the routing. The former consists in assigning the nodes to a specific cluster (or vehicle). The latter, once the assignment is performed, consists in determining, for each cluster, the Hamiltonian cycle defined on the subgraph that contains the cluster nodes. The idea underlying the Markov Decision Process proposed to model the problem is the following:  starting from a feasible solution, the agent picks a node and assigns it to another feasible cluster. The reward the agent receives is the difference between the objective function in two successive assignments. In other words, the goal of the model is for the agent to learn just an assignment strategy. The routing part, instead, is performed by a traditional heuristic procedure, which makes a good estimate of the objective function.

Using the notation introduced in the paragraph 2, let $n$ be the number of nodes, let $m$ be the number of vehicles (clusters), let $l$ be the number of features. We define $A\in\mathbb{R}$ $^{n\times n}$ the adjacency matrix,$X\in \{0,1\} ^{l\times n}$the assignment matrix (which is a boolean matrix where the element $x_{i_v}=1$ if and only if the node $i$ is assigned to the cluster $v$). In addition, let $D\in\mathbb{R}$ $^{l\times n}$ the demand matrix ,$Q\in \mathbb{R}$ $^{l\times m}$ the capacity matrix. In conclusion $y\in\mathbb{R}$ $^m$ is the vector containing the approximated cost of each cluster. When we mention the cost of a cluster, we refer to the cost of the optimal solution of a Travelling Salesman Problem (TSP), having as nodes those of that specific cluster. However, as already mentioned in the paragraph 2, finding the cost of a TSP is a NP-Hard problem. For this reason, $y_v$, which is the approximation of the cost of the cluster $v$, is computed using a traditional approximate heuristic (Christofides \cite{christofides1976worst}), which guarantees a good upper bound (not worse than $\frac{3}{2}$ the optimal solution) computed efficiently. The role of the $y$ variable, in fact, is to provide an estimate of the quality of a given cluster, and so it is not necessary to solve the TSP to the optimum.

After having introduced the basic notation, we can characterize the finite MDP for the VRP using the tuple $(\mathcal S, \mathcal A, \mathcal R, \mathcal P)$, 
where $\mathcal S$ is the set of states, $\mathcal A$ the set of actions, $\mathcal{R}:\mathcal{S}\times \emph{A}\rightarrow \mathbb{R}$ is the reward function and $\mathcal{P}:\mathcal{S}\times \emph{A}\rightarrow \mathcal{S}$ the transition function.
The state $s_t\in \mathcal{S}$  captures all the relevant information regarding the current iteration, in order to respect
the Markov Property and to have a fully observable system state. In addition, the state has the role of providing the agent all the relevant information about the problem, so that he can learn a strategy based on the peculiar structure of the problem. In this specific setting, the state $s_t$ is a tuple defined as $s_t=(A,X_t,y_t,D,Q)$. The adjacency matrix, as well as the demand and the capacity matrix, are not indexed by the iteration $t$ since their values are independent of the assignment chosen. Given a state $s_t$, the agent performs the action $a_t$, which is made up of two consecutive steps: first, the agent picks a node other than the depot, then it assigns the chosen node to a feasible cluster, meaning that the new assignment must satisfy the capacity constraints. We suppose the agent chooses to move the node $j$ to the cluster $k$. Then, the next state $s_{t+1}$ is given by the tuple $s_
{t+1}=(A,X_{t+1},y_{t+1},D,Q)$. As we mentioned before, $A,D,Q$ do not vary from one iteration to the next. Instead $x_{iv}^{t+1}$, which is a generic element of the assignment matrix is given by:
\begin{equation}
x_{iv}^{t+1}  = \left\{\begin{array}{cr}
      x_{iv}^{t} \text{ if } i \neq j\cap v\neq k \text{  or if  } x_{jk}^t=1 \\
      1 \text{   if         } i=j\cap v=k\\
      0 \text{   if         } i=j\cap v\neq k\
\end{array} \right.
\end{equation}
Similarly, ${y_v}^{t+1}$ is the approximated cost of the optimal solution of a TSP, having as nodes those of the cluster $v$, related to the assignment ${X}^{t+1}$.
In conclusion, the reward $r^{t+1}$ is obtained as the difference between the approximated cost of the VRP in the successive two states, which means:

\begin{equation}
r^{t+1} =\sum_{v=1}^{m} {y_v}^{t+1}-{y_v}^{t}
\end{equation}

\section{Agents deep models}\label{sec:Agents deep models}
In this section, we describe the neural network architecture behind the deep reinforcement learning algorithm proposed. Since the method proposed belongs to the actor-critic class, two networks need to be specified: the actor for estimating the policy $\pi_{\theta(s_t)}$ and the critic for estimating the state-value function $V_\omega(s_t)$. In this specific setting, the actor does not perform a single action, but two consequential actions: the choice of the node and then the choice of the cluster. For this reason, the neural architecture for estimating the policy $\pi_\theta(s_t)$ is composed by two different neural networks, that we have named respectively Actor 1 and Actor 2. Hence, Actor 1 estimates $\pi_{{\theta}_1(s_t)}$, which represents the probability of choosing every node given the state $s_t$. Once the generic node $i$ is chosen according to $\pi_{{\theta}_1(s_t)}$, Actor 2 estimates $\pi_{{\theta}_2 (s_t|node=i)}$, which is the probability of assigning every cluster to the node i. Regarding the weights of the two networks, we used the notation $\theta_1$ for Actor 1 and $\theta_2$ for Actor 2 to emphasize that the two deep neural networks do not share parameters. 
The neural architecture proposed has the aim of being flexible to the size of the VRP instance. In fact, according to the state defined in the paragraph 3, the size of the matrices $(A,X_t,y_t,D,Q)$  varies among different instances since every problem differs in the number of nodes, clusters and features. This difficulty is typical of problems whose inputs can be represented as graphs. In fact, graphs nodes, in contrast to pixel images, do not have a fixed number of neighboring units nor the spatial order among them is fixed \cite{zhang2019graph}. For this reason, in recent years, many architecture belonging to geometric deep leaning \cite{bronstein2017geometric} and graph convolution neural networks (GCNN)\cite{zhang2019graph} have been proposed to tackle the peculiarity of the graph structure. However,the target problem, as mentioned in the paragraph 2, is characterized by having a complete graph $G$. As a consequence, the adjacency matrix $A$ is dense and its structure is similar to a pixel image, if $A$ is appropriately scaled. In this context, a classical Convolutional Neural Network(CNN) with a flexible padding dimension can be used, and is indeed the approach followed in this paper.

\begin{figure}[h]
\centering
   \includegraphics[scale=0.5]{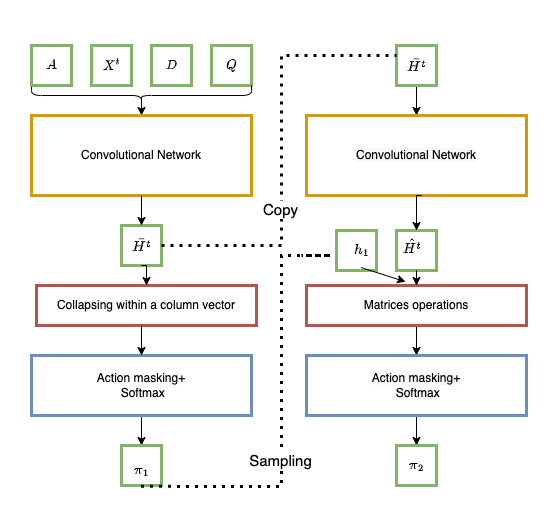}
    \caption{The Actor Model}
    \label{fig: The actor model}
\end{figure}

\begin{figure}[h]
    \centering
    \includegraphics[scale=0.6]{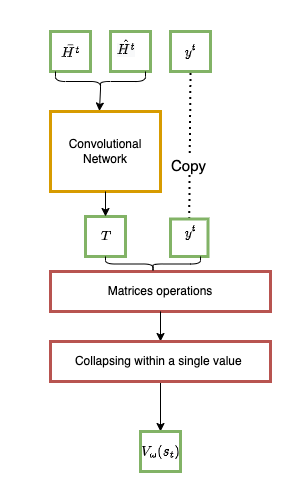}
    \caption{The Critic Model}
    \label{fig: The critic model}
\end{figure}

The Actor 1 model is composed by a deep Convolutional neural network (CNN), followed by a Softmax function. The CNN takes as input all the elements of the state $s_t$ apart from the vector $y_t$ and it produces as output the embedded matrix $\Bar{H^t} \in \mathbb{R}$ $^{m\times n}$, which represents the encoding of all the relevant information about the problem in a given iteration $t$. Then, the matrix $\Bar{H^t}$ is collapsed within a column vector of dimension $n$ and the Softmax function is applied to obtain $\pi_{{\theta}_1 (s_t)}$. Afterwards, a node is chosen according to $\pi_{{\theta}_1(s_t)}$. We introduce the  incidence vector ${h_1}^t \in\{0,1\} ^{n}$, whose elements are all $0$, apart from the one relative to the node chosen. At this stage the Actor 2, which is composed by a CNN, followed by a Softmax function as well, takes as input the couple $(h_1^t,\Bar{H^t})$. The CNN produces as output another embedded matrix $\Hat{H^t} \in \mathbb{R}$ $^{m\times n}$.Then, similarly to Actor 1, the matrix is collapsed within a columns vector of dimension $m$ and the Softmax function is applied to obtain $\pi_{{\theta}_2(s_t)}$. At this regard, we point out that, before applying the Softmax function, a mask is applied to force Actor 2 to choose a  feasible cluster. A graphic abstract representation of the Actor model is reported in figure \ref{fig: The actor model}. Lastly, the Critic, figure \ref{fig: The critic model}, is composed by a deep CNN as well. It takes as input the tuple $(y^t,\Bar{H^t},\Hat{H^t})$, it performs convolution operations, as well as simple matrix operations. The output of these operations is a matrix, that is collapsed within a scalar output through simple sum operations. The complete representation of the neural networks are reported in the appendix.

\section{RL Algorithm} \label{sec:algorithms}
In this section we describe in detail the reinforcement learning algorithm proposed. As already mentioned in the paragraph 2, it is an adaptation of the Proximal Policy Optimization (PPO) algorithm, more specifically the Adaptive KL Penalty Coefficient version.

The proposed algorithm is of the improvement type, meaning that the agent has the goal of improving a starting feasible solution. Therefore, the task the agent has to perform does not have a  terminal state since improving a solution is a process that lasts until the optimum is found, and it is impossible given a solution to verify whether it is really the optimum. For this reason, the task  belongs to the so called non episodic or continuous task. However,we can turn it  into an episodic task if we consider an episode as a fixed small number of steps $T$ in which the agent has to take actions to improve his current state. In this way, each episode does have a terminal state, which corresponds to the state $s_T$. Differently from the "real" episodic tasks, in our case the terminal state does not have a specific meaning within the problem. Therefore, also the reward $r_T$ does not have a particular role if compared to the rewards to the previous steps. For this decision, the choice we made is to store, for each roll-out, only the advantage function at the first step $\hat A_1$. A detailed description of the algorithm, in term of pseudo-code is reported below (Algorithm \ref{ref:PPO_Algorithm_VRP}).

\begin{algorithm}[H]\caption{Proximal Policy Optimization for VRP}\label{ref:PPO_Algorithm_VRP}
\DontPrintSemicolon
\SetAlgoLined
Choose the hyperparameters ${\text{learning rate}}_{\text{actor}},  {\text{learning rate}}_{\text{critic}},{epochs}_{actor}, {epochs}_{critic}$\;
 Initialize $\theta=\theta_0$, $\omega=\omega_0$, $\beta=\beta_0$\; Given a feasible initial solution compute $s_0=(A,X_0,y_0,D,Q)$\;
 Initialize two empty buffer $B_{iteration}, B_{rollout}$
 {

 \For{iteration =1, \dots, $k$}{
 \For{rollout = 1, \dots, N}{
 Clear buffer $B_{rollout}$\;
 \For{steps = 1, \dots, T}{
  Run $\pi_{{\theta_{1_t}} (s_t)}$ \;
  Choose the node according to $\pi_{{\theta_1} (s_t)}$: so compute the incidence vector $h_{1_t}$\;
  Retrieve the matrix $\bar{H_t}$\;
  Run $\pi_{\theta_{2_t | {h_{1_t}}}}(s_t)$\;
  Choose the cluster according to $\pi_{{\theta_2}(s_t)}$: so compute  $h_{2_t}$\;
  Compute $\pi(s_t)=\pi_{\theta_{1_t}}(s_t)\cdot \pi_{\theta_{2_t} | h_{1_t}}(s_t) $\;
  Retrieve the matrix $\hat{H_t}$\;
  Compute the state-value Function $V_{\omega_t} (s_t)$\;
  Compute the reward $r_t$}
  Compute the advantage function estimate $\hat A_1= \sum_{t=1}^{T} \gamma^t r_{t} -V_{\omega_1} (s_1)$ \;
  Store $\pi(s_1)$, $\bar{H_1}$, $\hat{H_1}$, $\hat A_1$ in the buffer $B_{rollout}$}
 
  Store the content of $B_{rollout}$ in $B_{iteration}$

 {
 Update the Actor parameter $\theta_1$ and $\theta_2$ according to (\ref{Actor objective function}), using data retrieved from $B_{iteration}$  \;
 Update the Critic parameter $\omega$  according to (\ref{Critic objective function}) using data retrieved from $B_{iteration}$  \;
 Clear buffer $B_{iteration}$\;

 }
 
 }

 }
\end{algorithm}

\section{Computational Experience}\label{sec:computational_experience}
In this section we describe the computational experiments performed on VRP instances. The main goals of the experiments are showing that the algorithm proposed has learning capabilities (i.e.the agent can learn strategies to decrease the objective function) and that the model, trained on a training set, has generalization capabilities on the test set. In addition, we want to show the  comparison between the quality of the performance obtained and the ones concerning the state of the art solver for Combinatorial optimization heuristics, which is $OR-TOOLS$ (https://developers.google.com/optimization).
The numerical tests were carried out in two stages.

\begin{figure}
\centering 
\begin{subfigure}{0.45\textwidth}
\centering 
  \includegraphics[width=\linewidth]{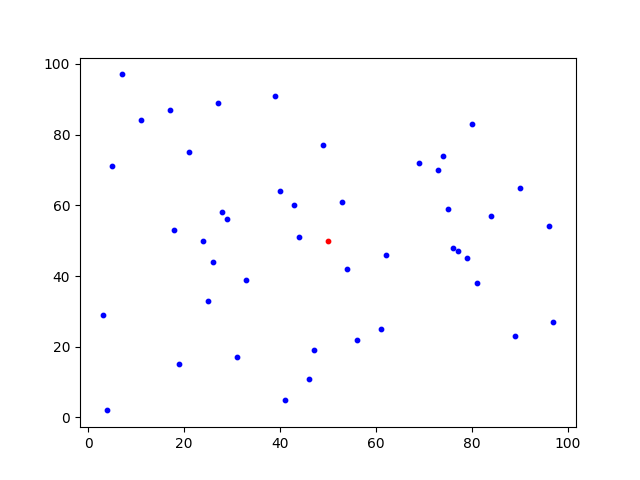}
  \caption{Uniform configuration}
  \label{fig:uniform}
\end{subfigure}%
\hfill
\begin{subfigure}{0.45\textwidth}
  \centering 
  \includegraphics[width=\linewidth]{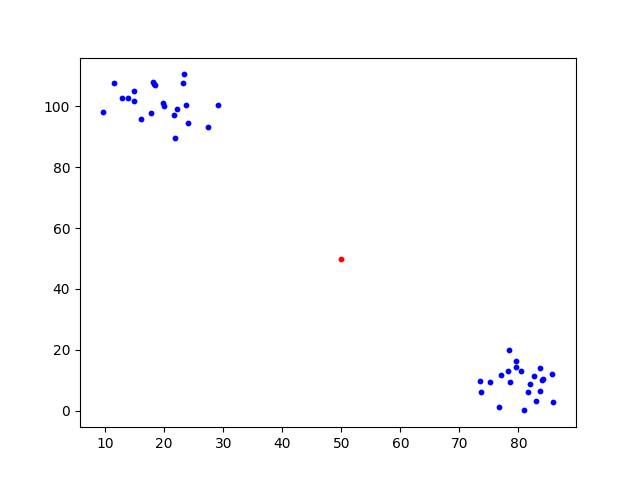}
  \caption{Two clusters configuration}
  \label{fig:clusters}
\end{subfigure}
\caption{Different types of nodes configuration.}
\label{fig:configurations_example}
\end{figure}

\textbf{-}During the first stage, we generated training and test instances belonging to  three classes with different characteristics to evaluate the behavior of the algorithm in specific situations. During this stage, for each class of instances, we compared the performances of the algorithm proposed with the $OR- TOOLS$ Methaheuristic TABOO SEARCH in terms of value of the objective function returned. We point out that the comparison was performed for different problem dimensions. Hence, for each class of instances and for every dimension, we computed the average of the objective function returned for $10$ instances. In conclusion, the run time of the two algorithms was not a matter of comparison; so we conducted the tests using the same run time for the two algorithms. 
The three classes of instances are the following:
\begin{itemize}
    \item[\textbf{C1}] In this class there is only one feature and the capacity of each vehicle related to that feature is $\infty$. In addition, we suppose the nodes are uniformly distributed in the two-dimensional space; more specifically we suppose $x_i \sim \mathcal{U}_{[0,100]\times[0,100]}$. The depot is located at the centre, in the position of coordinates $(50;50)$ (Figure \ref{fig:uniform}). The reason underlying this class of instances is to show the performances of the algorithm without the complexity of the capacity constraints.
    
    \item[\textbf{C2}] In this class there is only one feature, the capacity of each vehicle is the same and the nodes are arranged in two separate clusters , with the depot in the middle between the two (Figure \ref{fig:clusters}).The reason underlying this class of instances is to verify if the algorithm is able to recognise a specific pattern and identify a proper assigning strategy.
    
    \item[\textbf{C3}] In this class there is only one feature, the capacity of each vehicle is the same and the nodes are uniformly distributed in the  two-dimensional space,i.e. $x_i \sim \mathcal{U}_{[0,100]\times[0,100]}$. The depot is located at the centre, in the position of coordinates $(50;50)$ (Figure \ref{fig:uniform}).The reason underlying this class of instances is to verify if the algorithm is able to identify a proper assigning strategy in a more generic setting.

\end{itemize}

\begin{figure}
\centering 
\begin{subfigure}{0.45\textwidth}
\centering 
  \includegraphics[width=\linewidth]{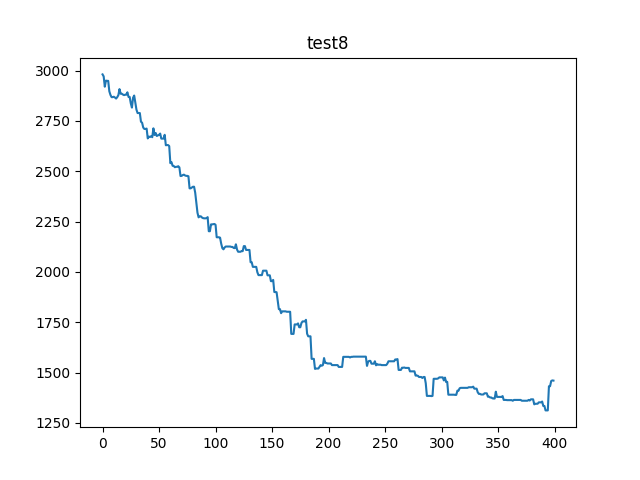}
  \caption{Objective function decrease in the case \textbf{C1} }
  \label{fig:uniform}
\end{subfigure}%
\hfill
\begin{subfigure}{0.45\textwidth}
  \centering 
  \includegraphics[width=\linewidth]{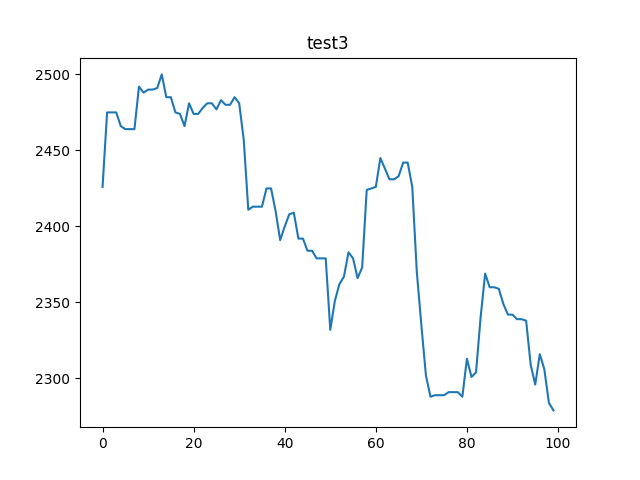}
  \caption{Objective function decrease in the case \textbf{C3}}
  \label{fig:objective function decrease}
\end{subfigure}
\caption{Decrease of the objective function for different configurations.}
\label{fig:Decreases of objective function}
\end{figure}

As shown in the table 1, the class of instances in which the algorithm performs the best is the class $\textbf{C1}$. In fact, for this class of instances,it provides a solution, which is only about $20-25\%$ worse than the one provided by the solver $OR-TOOLS$ with the same run time.
Regarding the other two class of instances, the gap between the algorithms grows in favor of $OR-TOOLS$ solver. However, we want to show that, even if the solution returned is much worse, the algorithm does learn a strategy (i.e. the objective function decreases, see figure \ref{fig:Decreases of objective function}).

\begin{table}[h]
\begin{subtable}{1\textwidth}
    \centering
    \begin{tabular}[width=.7\linewidth]{c|c|c|c|c}
    \textbf{Number of nodes} &\textbf{Number of vehicles}& \textbf{Run time}& \textbf{OR-TOOLS}& \textbf{Reinforcement learning} \\
    \hline
    Between 30 and 40 & Between 3 and 4 & 30 seconds & 501.30 & 629.20 \\
    Between 60 and 70 & Between 5 and 6 & 60 seconds & 665.0 & 821.80 \\

     & & \\
    \end{tabular}
    \caption{Performance on the class of instances \textbf{C1}}
    \label{tab:C1 performances}
\end{subtable}

\bigskip

\begin{subtable}{1\textwidth}
    \centering
    \begin{tabular}[width=.7\linewidth]{c|c|c|c|c}
    \textbf{Number of nodes} &\textbf{Number of vehicles}& \textbf{Run time}& \textbf{OR-TOOLS}& \textbf{Reinforcement learning} \\
    \hline
    Between 30 and 40 & Between 3 and 4 & 30 seconds & 313.90 & 878.60  \\
    Between 60 and 70 & Between 5 and 6 & 60 seconds & 359.00 & 1410.69 \\

     & & \\
    \end{tabular}
    \caption{Performance on the class of instances \textbf{C3}}
    \label{tab:C2 performances}
\end{subtable}

\bigskip

\begin{subtable}{1\textwidth}
    \centering
    \begin{tabular}[width=.7\linewidth]{c|c|c|c|c}
    \textbf{Number of nodes} &\textbf{Number of vehicles}& \textbf{Run time}& \textbf{OR-TOOLS}& \textbf{Reinforcement learning} \\
    \hline
    Between 30 and 40 & Between 3 and 4 & 30 seconds & 506.20 & 907.30 \\
    Between 60 and 70 & Between 5 and 6 & 60 seconds & 719.60 & 1637.40 \\

     & & \\
    \end{tabular}
    \caption{Performance on the class of instances \textbf{C3}}
    \label{tab:C3 performances}
\end{subtable}
\caption{Comparison of the performance with $OR-TOOLS$ solver}
\end{table}

\begin{figure}[h]
    \centering
    \includegraphics[scale=0.5]{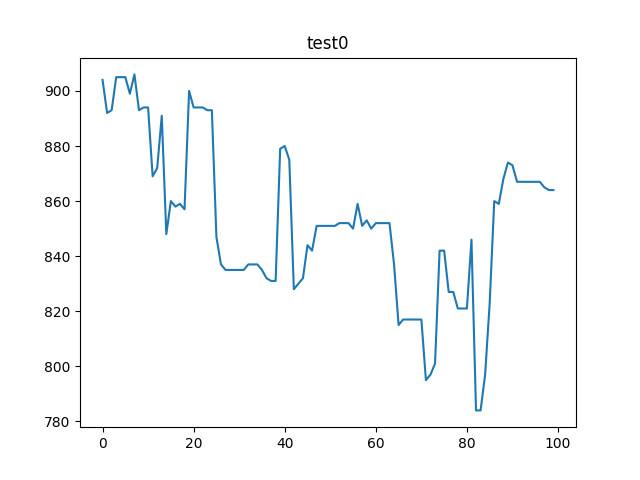}
    \caption{Objective function decrease in the case of $CVRPLIB$ instances}
    \label{fig: CVRPLIB}
\end{figure}

\textbf{-}During the second stage, we performed the training and test process on another class of instances, retrieved from the repository $CVRPLIB$, which is the state of the art repository for the Capacitated vehicle routing problem. In fact, it contains instances regarding real distribution problems. Among all the instances present there, we picked 85 of them for the training and 15 for the test. During this stage the goal of the experiments is not to compare the performance of the algorithm proposed with the $OR-TOOLS$ solver, but rather to show that the model proposed has learning capabilities even with real instances. At this regard, we show a figure(See \ref{fig: CVRPLIB}) regarding the behavior of the algorithm on the test set, in terms of value of the objective function throughout the iterations.

In conclusion, the computational experience shows that the algorithm does learn strategies to improve the objective function, but the results produced are still far from the ones produced by the solver $OR-TOOLS$. However,there are future research perspectives for the work proposed. These mainly concern the use of graph neural networks instead of convolutional neural networks and the possibility of the inclusion of a neural network for estimating the cost of each TSP.

\section{Conclusions}
\label{sec:Conclusions}
In this paper we investigated how to solve the Capacitated Vehicle routing Problem through Deep Reinforcement learning. Although we applied a different approach compared to the traditional one, we were inspired by two classic principles of  VRP heuristics: the improving mechanism and the two stage method. At this regard, the aim of the work was to build an improving algorithm capable of learning the assignment part of the VRP for a wide class of instances.

We firstly formulated the CVRP as a Markov Decision Process (MDP), identifying for the specific problem the tuple $(\mathcal S, \mathcal A, \mathcal R, \mathcal P)$. After that, we included the MDP identified in a RL algorithm. We decided to use, as algorithms, an adaptation of the PPO method, which belong to the class of Actor-Critic.
In a second phase we addressed the problem of which typology of neural architecture to use for the Actor and Critic respectively. The neural networks selected were supposed to have at least two fundamental characteristics: as first, they were supposed to be flexible to the input dimension. In addition, they were supposed to take advantage of the particular structure of the input, namely that of a complete graph. For this reason, the neural architecture chosen is made of Convolutional deep neural networks with a flexible padding dimension.
In conclusion the last part is about computational results. the primary goal of this part was to show that  an improving deep reinforcement learning algorithm for the VRP is possible, i.e. the algorithm proposed has learning capabilities. In addition, we showed the comparison between the results obtained and the ones obtained through $OR-TOOLS$ solver. The comparison shows how $OR-TOOLS$ solver still guarantees better results in terms of solution returned. However, there are future research perspectives for the work proposed. These mainly concern the inclusion of a neural network for estimating the cost of each TSP, instead of using the Christofides algorithm and the use of graph neural networks instead of convolutional neural networks.

\clearpage
\bibliographystyle{apalike}
\bibliography{biblio}

\clearpage

\section*{Appendix}
\subsection*{Hyper-parameters setting}

\bigskip
\bigskip

\begin{center}
    \begin{tabular}{|c|c|}
    \hline
    Actor Optimizer & Adam \\
    \hline
    Actor Learning Rate & $1 \times 10^{-5}$\\
    \hline
    Critic Optimizer & Adam \\
    \hline
    Critic Learning Rate & $1 \times 10^{-5}$\\
    \hline
    Number of Kernels of CNN & 27 \\
    \hline
    Number of output channel of CNN & 4 \\
    \hline
    \end{tabular}
\end{center}

\bigskip 
\bigskip

\subsection*{Detailed neural architectures}
\begin{figure}[h]
    \centering
    \includegraphics[scale=0.4]{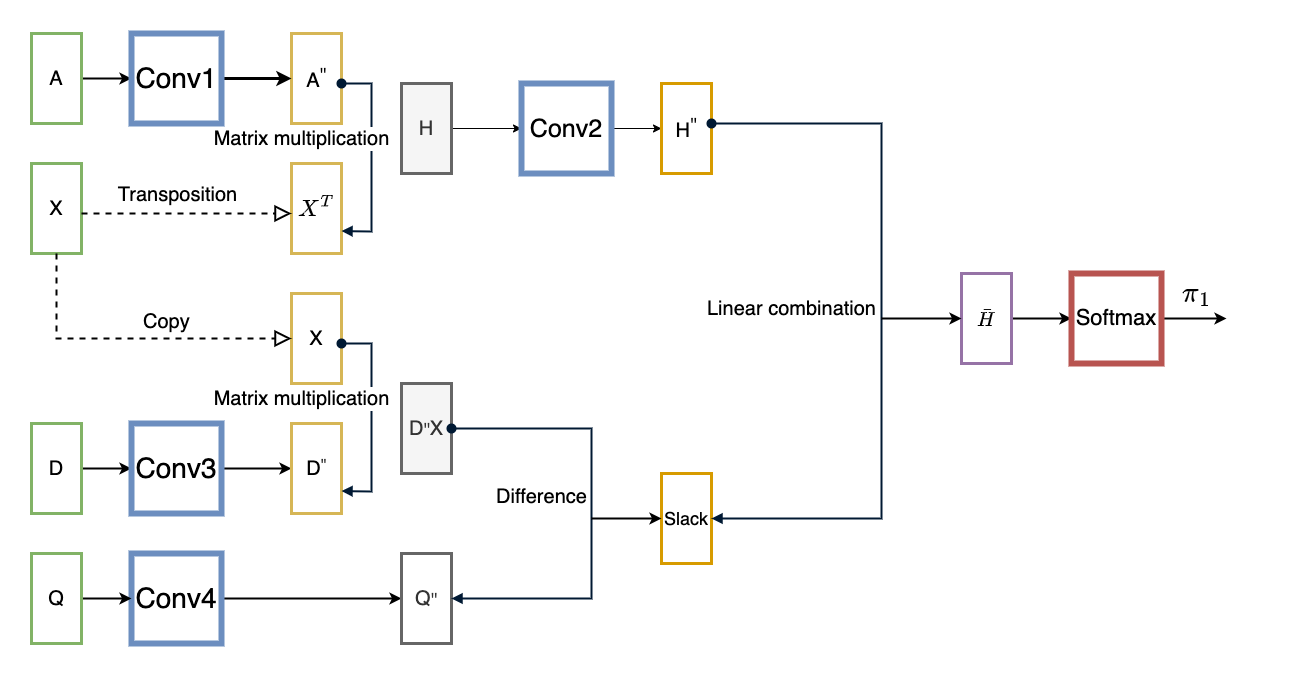}
    \caption{The Actor 1 model}
    \label{fig: The Actor 1 model}
\end{figure}

\begin{figure}[h]
    \centering
    \includegraphics[scale=0.4]{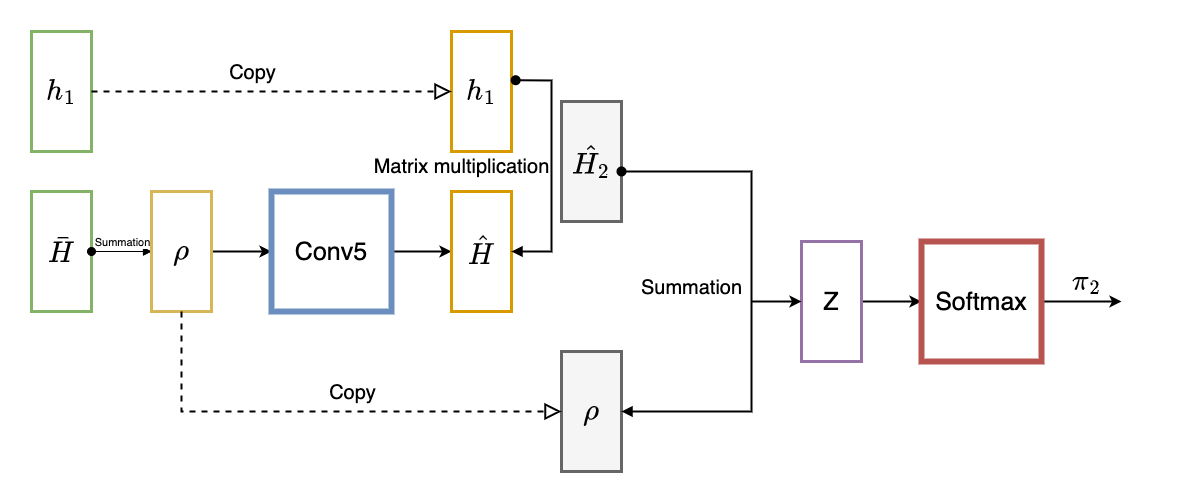}
    \caption{The Actor 2 Model}
    \label{fig: The actor 2 model}
\end{figure}

\begin{figure}[h]
    \centering
    \includegraphics[scale=0.35]{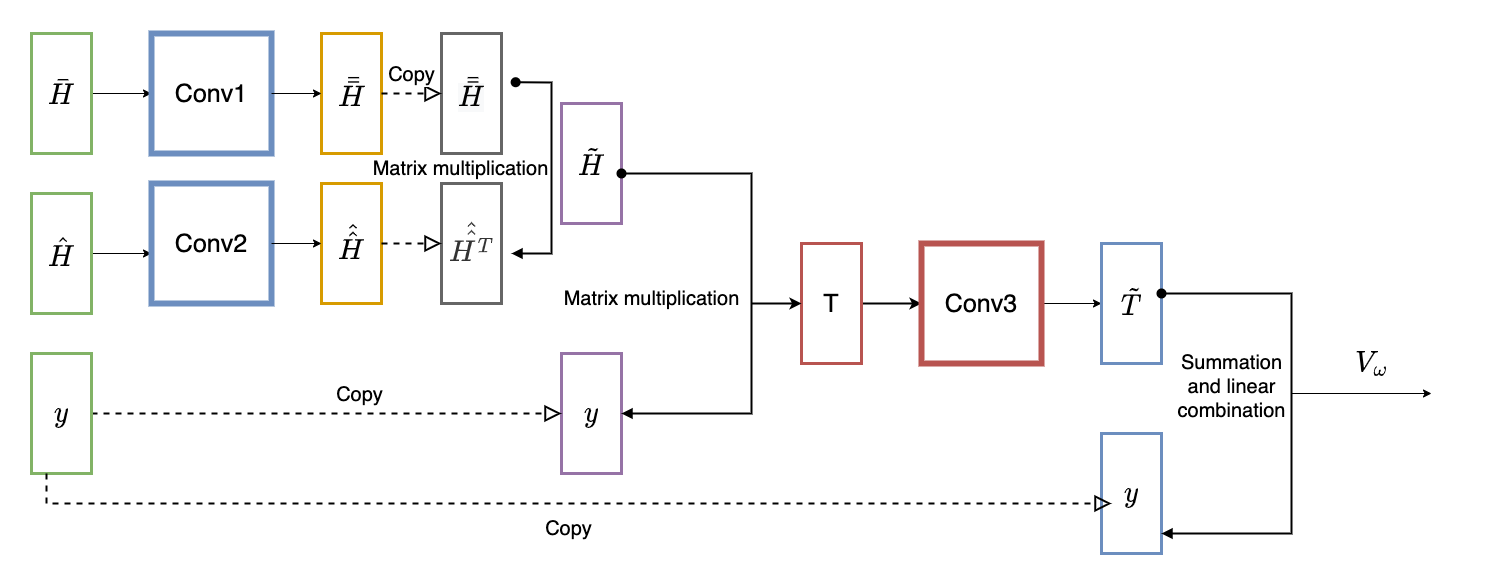}
    \caption{The Critic complete model}
    \label{fig: The Critic complete model}
\end{figure}

\end{document}